\documentclass{birkjour}

\usepackage[unicode]{hyperref}
\usepackage{latexsym}
\usepackage{indentfirst}
\usepackage{amsxtra}
\usepackage{amssymb}
\usepackage{amsthm}
\usepackage{amsmath}
\usepackage{mathrsfs} 
\usepackage[capitalise]{cleveref}

\newcommand{\overbar}[1]{\mkern1.5mu\overline{\mkern-1.5mu#1\mkern-1.5mu}\mkern 1.5mu}

\newtheorem{theor}{Theorem}
\newtheorem*{theor*}{Theorem}
\newtheorem{prop}[theor]{Proposition}
\newtheorem{lemma}[theor]{Lemma}
\newtheorem{cor}[theor]{Corollary}
\newtheorem*{cor*}{Corollary}
\theoremstyle{definition}               %stile roman 
\newtheorem{defin}[theor]{Definition}

\DeclareMathOperator{\Sym}{Sym}
\DeclareMathOperator{\Aut}{Aut}
\DeclareMathOperator{\End}{End}

\DeclareMathOperator{\id}{id}

\DeclareMathOperator{\lcm}{lcm}
\DeclareMathOperator{\ind}{i}
\DeclareMathOperator{\per}{p}

\newcommand{\lambdaa}[2]{\lambda_{#1}{#2}}
\newcommand{\lambdaam}[2]{\lambda^{-1}_{#1}{#2}}
\newcommand{\rhoo}[2]{\rho_{#1}{#2}}
\newcommand{\alphaa}[3]{\alpha^{#1}_{#2}{#3}}
\newcommand{\betaa}[3]{\beta^{#1}_{#2}{#3}}
\newcommand{\lambdai}[3]{\lambda^{#1}_{#2}{#3}}
\newcommand{\rhoi}[3]{\rho^{#1}_{#2}{#3}}
\newcommand{\indd}[1]{\ind{#1}}
\newcommand{\perr}[1]{\per{#1}}

\begin{document}
	
	\title[The matched product of the solutions of finite order]
	{The matched product of the solutions to the Yang-Baxter equation of finite order}

	%----------Author 1
	\author[F. Catino]{Francesco~CATINO}
	
	\address{%
		Department of Mathematics and Physics ``Ennio De Giorgi''\\ 
		University of Salento \\
		Via per Arnesano\\
		73100 Lecce (Le)\\
		Italy\\
		OrcID profile: \url{https://orcid.org/0000-0002-6683-1945}}
	
	\email{francesco.catino@unisalento.it}
	
	\thanks{This work was partially supported by the Department of Mathematics and Physics ``Ennio De Giorgi" -- University of Salento. The authors are members of GNSAGA (INdAM)}
	%----------Author 2
	\author[I. Colazzo]{Ilaria~COLAZZO}
	\address{%
		Department of Mathematics and Physics ``Ennio De Giorgi''\\ 
		University of Salento \\
		Via per Arnesano\\
		73100 Lecce (Le)\\
		Italy\\
		OrcID profile: \url{https://orcid.org/0000-0002-2713-0409}}
	\email{ilaria.colazzo@unisalento.it}
	%----------Author 3
	\author[P. Stefanelli]{Paola~STEFANELLI}
	\address{%
		Department of Mathematics and Physics ``Ennio De Giorgi''\\ 
		University of Salento \\
		Via per Arnesano\\
		73100 Lecce (Le)\\
		Italy
		OrcID profile: \url{https://orcid.org/0000-0003-3899-3151}}
	\email{paola.stefanelli@unisalento.it}
	\subjclass{Primary 16T25; Secondary 81R50, 16Y99, 16N20}
	
	\keywords{Quantum Yang-Baxter equation, set-theoretical solution, brace, semi-brace}
	
	\date{\today}
%	\dedicatory{To my boss}
	%%% ----------------------------------------------------------------------
	
	\begin{abstract}
			In this work, we focus on the set-theoretical solutions of the Yang-Baxter equation which are of finite order and not necessarily bijective. We use the matched product of solutions as a unifying tool for treating these solutions of finite order, that also include involutive and idempotent solutions. In particular, we prove that the matched product of two solutions $r_S$ and $r_T$ is of finite order if and only if $r_S$ and $r_T$ are. Furthermore, we show that with sufficient information on $r_S$ and $r_T$ we can precisely establish the order of the matched product. Finally, we prove that if $B$ is a finite semi-brace, then the associated solution $r$ satisfies $r^n=r$, for an integer $n$ closely linked with $B$.
	\end{abstract}
%%% ----------------------------------------------------------------------
\maketitle
%%% ----------------------------------------------------------------------	
\section{Introduction}
	The Yang-Baxter equation is a fundamental tool in several different fields of research such as statistical mechanics, quantum group theory, and low-dimensional topology. Named after the authors of the first papers in which the equation arose, Yang \cite{Ya67} and Baxter \cite{Ba72}, its study has been an extensive research area for the past sixty years.
    In 1992, V. Drinfel$'$d \cite{Dr90} suggested focusing on a specific class of solutions: the \emph{set-theoretical solutions} or \emph{braided sets}. Namely, given a set $X$, a set-theoretical solution, shortly a solution, is a map $r: X\times X \to X \times X$ such that the following condition 
    \begin{align*}
        \left(r\times \id_X\right)\left(\id_X\times r\right)\left(r\times \id_X\right) = \left(\id_X\times r\right)\left(r\times \id_X\right)\left(\id_X\times r\right)
    \end{align*}
    is satisfied. If $r$ is such a solution on X, for $x, y \in X$, define the maps ${\lambda_x : X \to X}$ and $\rho_y : X \to X$ by $r\left(x, y\right) = \left(\lambdaa{x}{(y)}, \rhoo{y}{(x)}\right)$. A solution $r$ is said to be \emph{left} (resp. \emph{right}) \emph{non-degenerate} if $\lambdaa{x}{}$ (resp. $\rhoo{x}{}$) is bijective, for each $x \in X$. The seminal papers by Etingof, Schedler, and Soloviev \cite{ESS99}, and Gateva-Ivanova and Van den Bergh \cite{GaB98}, laid the foundations for studying the class of non-degenerate solutions that are also involutive. A solution $r$ on $X$ is said to be involutive if $r^2 = \id$. Such solutions have been intensively studied, see \cite{GaMa08,CeJR10,Ve16,BaCJO17,CaCP18} just to name a few. In particular, Rump in \cite{Ru07a} introduced \emph{braces}, ring-like structures, for studying involutive non-degenerate solutions. As reformulated by Ced\'o, Jespers, and Ok\'ninski in \cite{CeJO14}, a left brace is a set $B$ with two operations $+$ and $\circ$ such that $(B,+)$ is an abelian group, $(B,\circ)$ is a group and the relation
    $a \circ\left(b+c\right) +a = a \circ b + a \circ c$
    holds for all $a,b,c \in B$. Braces have been widely studied, see for instance \cite{Ru07b, CCoSt16, Ga18, Sm18, BaCeJ18, Ba18, CeGaSm18, CeSmVe19}.
    Soloviev in \cite{So00} and Lu, Yan, and Zhu in \cite{LuYZ00} studied bijective not necessarily involutive solutions. Such solutions have been relatively investigated \cite{Ya16,Ya18,ElNeTs19} and have applications in knot theory, see \cite{Ne11} and the references therein. Guarnieri and Vendramin in \cite{GVe17} introduced \emph{skew braces}, structures, including braces, useful for studying this class of solutions. A set $B$ with two operations $+$ and $\circ$ is a skew left brace if $\left(B,+\right)$ and $\left(B,\circ\right)$ are groups and the condition
    $a \circ\left(b+c\right) = a \circ b - a + a \circ c$
    is satisfied for all $a,b,c\in B$.
    Further advancements in the field of skew braces relating to Hopf-Galois structures can be found in \cite{SmVe18, Ch18, NZ19}, whereas \cite{CCoSt19}, an extension of \cite{CCoSt15}, partially answered the extension problem in a simplified case.
    It is worth mentioning that in literature bijective solutions are usually defined on finite sets. Under this assumption, for each bijective solution $r$, there exists an integer $n$ such that $r^n=\id$.
    In \cite{Le17}, Lebed drew attention on idempotent solutions that, although of little interest in physics, provide a tool for dealing with very different algebraic structures ranging from and free (commutative) monoids to factorizable monoids, and from distributive lattices to Young tableaux and plactic monoids. Namely, given a set $X$, a solution $r$ is said to be idempotent if $r^2=r$.  The question arises whether there is an algebraic structure similar to the brace structure useful for studying solutions not necessarily bijective. In \cite{CCoSt17}, we gave an initial answer to the question by introducing \emph{semi-braces}.  A left (cancellative) semi-brace is a set $B$ with two operations $+$ and $\circ$ such that $\left(B,+\right)$ is a left cancellative semigroup, $\left(B,\circ\right)$ is a group and $a\circ\left(b+c\right) = a \circ b + a \circ\left(a^-+c\right)$ holds for all $a,b,c\in B$, where $a^-$ denotes the inverse of $a$ in $\left(B,\circ\right)$. 
    Later, in \cite{VAJe19}, this algebraic structure has been generalized by weakening the hypotheses of the additive semigroup, in particular by removing the assumption of left cancellativity.
    The connection between involutive, bijective and idempotent solutions, at least on a finite set, is that for every solution $r$ there exist two non-negative integers $i$ and $p$ such that $r^{p+i} = r^{i}$. Such a solution $r$ is said to be of \emph{finite order}, and the minimal non-negative integers that satisfy such relation are said to be \emph{index} and \emph{period} and are denoted by $\indd{\left(r\right)}$ and $\perr{\left(r\right)}$, respectively.
    
    The aim of this work is to show how the matched product of solutions is a unifying tool for treating solutions of finite order.
    Let $r_S$ and $r_T$ be solutions on the sets $S$ and $T$, respectively. In \cite{CCoSt18x} we define a new solution on the cartesian product of $S$ and $T$ called the matched product of $r_S$ and $r_T$ denoted by $r_S\bowtie r_T$.
    Specifically, we prove that $r_S\bowtie r_T$ is a solution of finite order if and only if $r_S$ and $r_T$ are. In particular, we can control the index and the period of the matched product of solutions $r_S$ and $r_T$ given their indexes and periods: 
    \medskip
    
    \noindent\textbf{Main result.} \emph{If $r_S$ and $r_T$ have fixed indexes and periods, be they $\indd{\left(r_S\right)}$, $\indd{\left(r_T\right)}$, $\perr{\left(r_S\right)}$, and $\perr{\left(r_T\right)}$ respectively, then the matched product of $r_S$ and $r_T$ has as the index the maximum of  $\indd{\left(r_S\right)}$ and $\indd{\left(r_T\right)}$ and as the period the least common multiple of periods $\perr{\left(r_S\right)}$ and $\perr{\left(r_T\right)}$. }
    
    \medskip
    \noindent As a consequence of this result we obtain that
    $r_{S}^{l} = \id$ and $r_{T}^{m} = \id$ if and only if $\left(r_{S}\bowtie r_{T}\right)^{n} = \id$,
    and $r_{S}^{l} = r_{S}$ and $r_{T}^{m} = r_{T}$ if and only if $\left(r_{S}\bowtie r_{T}\right)^{n} = r_{S}\bowtie r_{T}$,
    for some integers $l, m, n$.
    This corollary includes the particular case of involutive ($r^2=\id$) and idempotent solutions ($r^2=r$) already provided in \cite[Corollary 5]{CCoSt18x}.
    The main result also sheds new light on solutions associated with semi-braces: we prove that for each solution $r_B$ associated with a finite semi-brace $B$ (under some minimal assumptions) $r_B^{n}=r_{B}$, for a certain $n \in \mathbb{N}$. Indeed in the first section, we extend to the general case the matched product $S\bowtie T$ of two semi-braces $S$ and $T$, which we have introduced in \cite{CCoSt18x} for cancellative semi-braces. We also show that the solution $r_{S\bowtie T}$ associated with the matched product of two semi-braces $S$ and $T$ is the matched product of the solutions $r_S$ and $r_T$. Since a semi-brace $B$ is the matched product of two ``trivial'' semi-braces and a skew brace $G$, we obtain that $r_B^{n}=r_{B}$, for an integer $n$.  Based on this result we show that $r_{B}^{3} = r_{B}$ if and only if the sub-skew brace $G$ is a brace, improving \cite[Theorem 5.1]{VAJe19}.

\section{Definitions and preliminary results}
	
	In \cite{CCoSt18x}, we introduced a new construction technique for solutions of the Yang-Baxter equation that allows one to obtain new solutions on the cartesian product of sets, starting from completely arbitrary solutions.\\
	Given a solution $r_S$ on a set $S$ and a solution $r_T$ on a set $T$, if $\alpha: T \to \Sym\left(S\right)$ and $\beta: S \to \Sym\left(T\right)$ are maps, set $\alpha_u:=\alpha\left(u\right)$, for every $u \in T$, and $\beta_a:=\beta\left(a\right)$, for every $a\in S$, then the quadruple
	$\left(r_S,r_T, \alpha,\beta\right)$ is said to be a \emph{matched product system of solutions} if the following conditions hold
	
	{\scriptsize
		\begin{center}
			\begin{minipage}[b]{.5\textwidth}
				\vspace{-\baselineskip}
				\begin{align}\label{eq:primo}\tag{s1}
					\alpha_u\alpha_v = \alpha_{\lambda_u\left(v\right)}\alpha_{\rho_{v}\left(u\right)}
				\end{align}
			\end{minipage}%
			\hfill\hfill\hfill
			\begin{minipage}[b]{.5\textwidth}
				\vspace{-\baselineskip}
				\begin{align}\label{eq:secondo}\tag{s2}
					\beta_a\beta_b=\beta_{\lambda_a\left(b\right)}\beta_{\rho_b\left(a\right)}
				\end{align}
			\end{minipage}
		\end{center}
		\begin{center}
			\begin{minipage}[b]{.5\textwidth}
				\vspace{-\baselineskip}
				\begin{align}\label{eq:quinto}\tag{s3}
					\rho_{\alpha^{-1}_u\!\left(b\right)}\alpha^{-1}_{\beta_a\left(u\right)}\left(a\right) = \alpha^{-1}_{\beta_{\rho_b\left(a\right)}\beta^{-1}_b\left(u\right)}\rho_b\left(a\right)
				\end{align}
			\end{minipage}%
			\hfill\hfill
			\begin{minipage}[b]{.5\textwidth}
				\vspace{-\baselineskip}
				\begin{align}\label{eq:sesto}\tag{s4}
					\rho_{\beta^{-1}_a\!\left(v\right)}\beta^{-1}_{\alpha_u\left(a\right)}\left(u\right) = \beta^{-1}_{\alpha_{\rho_v\left(u\right)}\alpha^{-1}_v\left(a\right)}\rho_v\left(u\right)
				\end{align}
			\end{minipage}
		\end{center}
		\begin{center}
			\begin{minipage}[b]{.5\textwidth}
				\vspace{-\baselineskip}
				\begin{align}\label{eq:terzo}\tag{s5}
					\lambda_a\alphaa{}{\betaa{-1}{a}{\left(u\right)}}{}= \alphaa{}{u}{\lambdaa{\alphaa{-1}{u}{\left(a\right)}}{}}
				\end{align}
			\end{minipage}%
			\hfill\hfill
			\begin{minipage}[b]{.5\textwidth}
				\vspace{-\baselineskip}
				\begin{align}\label{eq:quarto}\tag{s6}
					\lambdaa{u}{\betaa{}{\alphaa{-1}{u}{\left(a\right)}}{}}=\betaa{}{a}{\lambdaa{\betaa{-1}{a}{\left(u\right)}}{}}
				\end{align}
			\end{minipage}
		\end{center}
	}
	\noindent for all $a,b \in S$ and $u,v \in T$.
	
	As shown in {\cite[Theorem 1]{CCoSt18x}}, any matched product system of solutions determines a new solution on the set $S\times T$. Specifically, if $\left(r_S,r_T, \alpha,\beta\right)$ is a matched product system of solutions, then the map $r:S{ \times} T\times S{ \times} T \to S{ \times} T\times S{ \times} T$ defined by
		\begin{align*}
			&r\left(\left(a, u\right), 
			\left(b, v\right)\right) := 
			\left(\left(\alphaa{}{u}{\lambdaa{\bar{a}}{\left(b\right)}},\, \betaa{}{a}{\lambdaa{\bar{u}}{\left(v\right)}}\right),\ \left(\alphaa{-1}{\overline{U}}{\rhoo{\alphaa{}{\bar{u}}{\left(b\right)}}{\left(a\right)}},\,  \betaa{-1}{\overline{A}}{\rhoo{\betaa{}{\bar{a}}{\left(v\right)}}{\left(u\right)}} \right) \right),
		\end{align*}
		where
		\begin{align*}
			\bar{a}:=\alphaa{-1}{u}{\left(a\right)}, \ &\bar{u}:= \betaa{-1}{a}{\left(u\right)},\\ A:=\alphaa{}{u}{\lambdaa{\bar{a}}{\left(b\right)}}, \ &U:=\betaa{}{a}{\lambdaa{\bar{u}}{\left(v\right)}}, \\ \overline{A}:=\alphaa{-1}{U}{\left(A\right)}, \ &\overline{U}:= \betaa{-1}{A}{\left(U\right)},
		\end{align*}
		for all $\left(a,u\right),\left(b,v\right)\in S\times T$, is a solution. This solution is called the \emph{matched product of the solutions} $r_S$ and $r_T$ (via $\alpha$ and $\beta$) and it is denoted by $r_S\bowtie r_T$.
		
	\noindent If $\left(r_{S}, r_{T}, \alpha, \beta\right)$ is a matched product system of solutions, we denote $\alphaa{-1}{u}{\left(a\right)}$ with $\bar{a}$ and $\betaa{-1}{a}{\left(u\right)}$ with $\bar{u}$, when the pair $\left(a,u\right) \in S \times T$ is
	clear from the context.

	There exists a special relation between the derived solution of $r_{S}\bowtie r_{T}$ and the derived solutions of $r_{S}$ and $r_{T}$.
	To show this relation, first recall that if $r$ is a left non-degenerate solution on a set $X$, then the map
	$r':X\times X\to X\times X$ defined by
	\begin{align*}
		r'\left(x,y\right):= \left(y, \lambdaa{y}{\rhoo{\lambdaam{x}{\left(y\right)}}{\left(x\right)}}\right),
	\end{align*}
	for all $x,y\in X$, is a solution called the \emph{derived solution} of $r$ (see for instance \cite{So00}).

	\begin{prop}
		Let $\left(r_{S}, r_{T}, \alpha, \beta\right)$ be a matched product system of solutions where $r_S$ and $r_T$ are left non-degenerate. Then it holds that
		\begin{align*}
			\left(r_{S}\bowtie r_{T}\right)' = 
			r_{S}'\times r_{T}'.
		\end{align*}
		\begin{proof}
			Note that 
			\begin{align*}
				\left(r_{S}\bowtie r_{T}\right)'\left(\left(a,u\right), \left(b,v\right)\right)
				&= \left(\lambda'_{\left(a,u\right)}\left(b,v\right), \rho'_{\left(b,v\right)}\left(a,u\right)\right)\\
				&= \left(\left(b,v\right), 
				\lambdaa{\left(b,v\right)}{\rhoo{\lambdaam{\left(a,u\right)}{\left(b,v\right)}}{\left(a,u\right)}}\right)\\
				&= \left(\left(b,v\right),
				\left(\lambdaam{b}{\rhoo{\lambdaam{a}{\left(b\right)}}{\left(a\right)}}, \lambdaam{v}{\rhoo{\lambdaam{u}{\left(v\right)}}{\left(u\right)}}\right)
				\right)
			\end{align*}
			since, set
			$\mathcal{A}:= \lambdaa{a}{\alphaa{}{\bar{u}}{\left(\alphaa{-1}{\bar{u}}{\lambdaam{a}{\left(b\right)}}\right)}}$
			and $\mathcal{U}:= \lambdaa{u}{\betaa{}{\bar{a}}{\left(\betaa{-1}{\bar{a}}{\lambdaam{u}{\left(v\right)}}\right)}}$, we have 
			\begin{align*}
				&\lambdaa{\left(b,v\right)}{\rhoo{\lambdaam{\left(a,u\right)}{\left(b,v\right)}}{\left(a,u\right)}}
				= 
				\lambdaa{\left(b,v\right)}{\rhoo{\left(\alphaa{-1}{\bar{u}}{\lambdaam{a}{\left(b\right)}}, \betaa{-1}{\bar{a}}{\lambdaam{u}{\left(v\right)}}\right)}{\left(a,u\right)}}\\
				&= \lambdaa{\left(b,v\right)}{\left(\alphaa{-1}{\overline{\mathcal{U}}}{\rhoo{\alphaa{}{\bar{u}}{\left(\alphaa{-1}{\bar{u}}{\lambdaam{a}{\left(b\right)}}\right)}}{\left(a\right)}}, 
					\betaa{-1}{\overline{\mathcal{A}}}{\rhoo{\betaa{}{\bar{a}}{\left(\betaa{-1}{\bar{a}}{\lambdaam{u}{\left(v\right)}}\right)}}{\left(u\right)}}\right)}\\
				&= \lambdaa{\left(b,v\right)}{\left(\alphaa{-1}{\overline{\mathcal{U}}}{\rhoo{\lambdaam{a}{\left(b\right)}}{\left(a\right)}}, 
					\betaa{-1}{\overline{\mathcal{A}}}{\rhoo{\lambdaam{u}{\left(v\right)}}{\left(u\right)}}\right)}\\
				&= \left(\lambdaa{b}{\alphaa{}{\bar{v}}{\left(\alphaa{-1}{\overline{\mathcal{U}}}{\rhoo{\lambdaam{a}{\left(b\right)}}{\left(a\right)}}\right)}},
				\lambdaa{v}{\betaa{}{\bar{b}}{\left(\betaa{-1}{\overline{\mathcal{A}}}{\rhoo{\lambdaam{u}{\left(v\right)}}{\left(u\right)}}\right)}}
				\right)\\
				&= \left(\lambdaa{b}{\alphaa{}{\bar{v}}{\alphaa{-1}{\bar{v}}{\rhoo{\lambdaam{a}{\left(b\right)}}{\left(a\right)}}}},
				\lambdaa{v}{\betaa{}{\bar{b}}{\betaa{-1}{\bar{b}}{\rhoo{\lambdaam{u}{\left(v\right)}}{\left(u\right)}}}}
				\right)\\
				&= \left(\lambdaa{b}{\rhoo{\lambdaam{a}{\left(b\right)}}{\left(a\right)}}, \lambdaa{v}{\rhoo{\lambdaam{u}{\left(v\right)}}{\left(u\right)}}\right),
			\end{align*}
			where the first equality holds since $\lambdaam{\left(a,u\right)}{\left(b,v\right)} = \left(\alphaa{-1}{\bar{u}}{\lambdaam{a}{\left(b\right)}}, \betaa{-1}{\bar{a}}{\lambdaam{u}{\left(v\right)}}\right)$ and  
			the second last equality holds because $\mathcal{A} = b$, $\mathcal{U} = v$, and so $\overline{\mathcal{A}} = \alphaa{-1}{\mathcal{U}}{\left(\mathcal{A}\right)} = \alphaa{-1}{v}{\left(b\right)} = \bar{b}$, and similarly $\overline{\mathcal{U}} = \bar{v}$.
		\end{proof}
	\end{prop}
	
	%%%%%%%%%%%%%%%%%%%%%%%%%%%%%%%%%%%%%%%%%%%%%%%%%%%%%%%%%%%%%%%%%%%%%%%
    In \cite{Ru07a}, Rump introduced braces  in order to obtain involutive solutions of the Yang-Baxter equation. With the aim of constructing bijective solutions, Guarnieri and Vendramin \cite{GVe17} introduced skew braces as a generalization of braces. In \cite{CCoSt17}, we introduced semi-braces that included skew braces and allowed to obtain left non-degenerate solutions. A semi-brace is a set $B$ with two operations $+$ and $\circ$ such that $\left(B,+\right)$ is a left cancellative semigroup, $\left(B,\circ\right)$ is a group, and 
    \begin{align}\label{eq:semibrace}
        a\circ\left(b+c\right) = a \circ b + a \circ \left(a^-+c\right)
    \end{align}
    holds for all $a,b,c\in B$ where $a^-$ is the inverse of $a$ in $\left(B,\circ\right)$.
    If $B$ is a semi-brace, in particular if it is a (skew) brace, we can define for all $a,b \in B$ 
    \begin{align*}
        \lambdaa{a}{\left(b\right)} := a \circ \left(a^-+b\right)\qquad \mbox{ and } \qquad \rhoo{b}{\left(a\right)}:=\left(a^-+b\right)^-\circ b,
    \end{align*}
    two maps $\lambda,\rho: B \to B^B$ such that $\lambda\left(a\right)\left(b\right):=\lambdaa{a}{\left(b\right)}$ and $\rho\left(a\right)\left(b\right):=\rho_a\left(b\right)$ and a map $r_B:B\times B \to B\times B$ defined by $r_B\left(a,b\right) :=\left(\lambdaa{a}{\left(b\right)},\rhoo{b}{\left(a\right)}\right)$. If $B$ is a brace or a skew brace or a semi-brace then $r_B$ is a solution, called \emph{solution associated} to $B$.
    Finally, the notion of semi-braces was generalized in \cite{VAJe19} by Jespers and Van Antwerpen.
    \begin{defin}[Definition 2.1 in \cite{VAJe19}]
        Let $B$ be a set with two operations $+$ and $\circ$ such that $\left(B,+\right)$ is a semigroup and $\left(B,\circ\right)$ is a group. One says that $\left(B, + , \circ \right)$ is \emph{left semi-brace} if
        \begin{align*}
            a \circ \left(b+c\right) = a \circ b + a\circ\left(a^- +c\right)
        \end{align*}
        for all $a, b, c \in B$. Here, $a^-$ denotes the inverse of $a$ in $\left(B, \circ \right)$. We call $\left(B, +\right)$ the \emph{additive semigroup} of the left semi-brace of $\left(B, + , \circ \right)$. If the semigroup $\left(B,+\right)$ has a pre-fix, pertaining to some property of the semigroup, we will also use this pre-fix with the left semi-brace.
    \end{defin}
    \noindent In particular following this definition, semi-braces previously introduced in \cite{CCoSt17} will be left cancellative semi-brace.
    
    We note that for a left semi-brace not necessarily left cancellative $B$ we can introduce $\lambda$ and $\rho$. Theorem $5.1$ in \cite{VAJe19} gives a sufficient condition for the map $r_B$, defined as in the left cancellative case, to be a solution.
    The following theorem states a necessary and sufficient condition to ensure that $r_B$ is a solution.

    \begin{theor}
        Let $\left(B,+,\circ\right)$ be a left semi-brace. The map $r_B:B\times B \times \to B\times B$, defined by $r_B\left(a,b\right):=\left(a\circ\left(a^-+b\right), \left(a^-+b\right)^-\circ b\right)$ for all $a,b \in B$, is a solution if and only if 
        \begin{align}\label{eq:condsolution}
            a + \lambdaa{b}{\left(c\right)}\circ\left(0 + \rhoo{c}{\left(b\right)}\right) = a + b\circ \left(0+c\right),
        \end{align}
        holds for all $a,b,c\in B$.
        \begin{proof}
            It is easily verified that $r_B$ is a solution if and only if 
            \begin{align}\label{eq:sollambda}
                &\lambdaa{a}{\lambdaa{b}{\left(c\right)}} = \lambdaa{\lambdaa{a}{\left(b\right)}}{\lambdaa{\rhoo{b}{\left(a\right)}}{\left(c\right)}}\\\label{eq:sollambdarho}
                &\lambdaa{\rhoo{\lambdaa{b}{\left(c\right)}}{\left(a\right)}}{\rhoo{c}{\left(b\right)}} =\rhoo{\lambdaa{\rhoo{b}{\left(a\right)}}{\left(c\right)}}{\lambdaa{a}{\left(b\right)}}\\\label{eq:solrho}
                &\rhoo{\rhoo{b}{\left(c\right)}}{\rhoo{\lambdaa{c}{\left(b\right)}}{\left(a\right)}} = \rhoo{c}{\rhoo{b}{\left(a\right)}}
            \end{align}
            for all $a,b,c \in B$. Let $a,b,c \in B$. 
            First recall that by Lemma 2.4 in \cite{VAJe19}
            \begin{align}\label{eq:lemma24}
                \forall a,b \in B \qquad a+b = a+0+b
            \end{align}
            and by Lemma 2.12 in \cite{VAJe19} $\lambda: B \to \End\left(B\right)$ is a homomorphism, i.e.,
            \begin{align}\label{eq:lambdahomomorphism}
                \lambdaa{a}{\lambdaa{b}{}} = \lambdaa{a\circ b}{},
            \end{align}
            for all $a,b \in B$.
            Furthermore note that 
            \begin{align}\label{eq:lambdacircrho}
                \lambdaa{a}{\left(b\right)}\circ \rhoo{b}{\left(a\right)} = a \circ b
            \end{align}
            for all $a,b\in B$ and then
            \begin{align*}
                \lambdaa{\lambdaa{a}{\left(b\right)}}{\lambdaa{\rhoo{b}{\left(a\right)}}{\left(c\right)}}
                =\lambdaa{\left(\lambdaa{a}{\left(b\right)}\right)\circ\left(\rhoo{b}{\left(a\right)}\right)}{\left(c\right)} = \lambdaa{a\circ b}{\left(c\right)} = \lambdaa{a}{\lambdaa{b}{\left(c\right)}},
            \end{align*}
            i.e., \eqref{eq:sollambda} holds. Moreover 
            \begin{align*}
                &\lambdaa{\rhoo{\lambdaa{b}{\left(c\right)}}{\left(a\right)}}{\rhoo{c}{\left(b\right)}} =
                \rhoo{\lambdaa{b}{\left(c\right)}}{\left(a\right)}\circ \left(\left(\rhoo{\lambdaa{b}{\left(c\right)}}{\left(a\right)}\right)^- +\rhoo{c}{\left(b\right)}\right)\\
                &=\left(a^-+\lambdaa{b}{\left(c\right)}\right)^- \circ \lambdaa{b}{\left(c\right)}\circ \left(\left(\lambdaa{b}{\left(c\right)}\right)^-\circ \left(a^-+\lambdaa{b}{\left(c\right)}\right)+\rhoo{c}{\left(b\right)}\right)\\
                &=\left(a^-+\lambdaa{b}{\left(c\right)}\right)^- \circ\left(a^-+\lambdaa{b}{\left(c\right)}+\lambdaa{b}{\left(c\right)}\circ\left(\left(\lambdaa{b}{\left(c\right)}\right)^-+\rhoo{c}{\left(b\right)}\right) \right)&\quad\mbox{by \eqref{eq:semibrace}}\\
                &=\left(a^-+\lambdaa{b}{\left(c\right)}\right)^- \circ\left(a^-+\lambdaa{b}{\left(c\right)}\circ\left(0+\rhoo{c}{\left(b\right)}\right)\right)&\quad\mbox{by \eqref{eq:semibrace}}\\
                &=\left(a^-+\lambdaa{b}{\left(c\right)}\right)^- \circ\left(a^-+b\circ\left(0+c\right)\right)&\quad\mbox{by \eqref{eq:condsolution}}
            \end{align*}
            and 
            \begin{align*}
                &\rhoo{\lambdaa{\rhoo{b}{\left(a\right)}}{\left(c\right)}}{\lambdaa{a}{\left(b\right)}} 
                =\left(\left(\lambdaa{a}{\left(b\right)}\right)^- + \lambdaa{\rhoo{b}{\left(a\right)}}{\left(c\right)}\right)^-\circ \lambdaa{\rhoo{b}{\left(a\right)}}{\left(c\right)}\\
                &=\left(\left(\lambdaa{a}{\left(b\right)}\right)^- +\rhoo{b}{\left(a\right)}\circ\left(\left(\rhoo{b}{\left(a\right)}\right)^-+c\right) \right)^-\circ\rhoo{b}{\left(a\right)}\circ\left(\left(\rhoo{b}{\left(a\right)}\right)^-+c\right)\\
                &=\left(\rhoo{b}{\left(a\right)}\circ\left(\rhoo{b}{\left(a\right)}\right)^-\circ\left(\lambdaa{a}{\left(b\right)}\right)^-+c\right)^-\circ\rhoo{b}{\left(a\right)}\circ\left(\left(\rhoo{b}{\left(a\right)}\right)^-+c\right)&\mbox{by \eqref{eq:semibrace}}\\
                &=\left(\left(\rhoo{b}{\left(a\right)}\right)^-\circ\left(\lambdaa{a}{\left(b\right)}\right)^-+c\right)^-\circ\left(\rhoo{b}{\left(a\right)}\right)^-\circ\rhoo{b}{\left(a\right)}\circ\left(\left(\rhoo{b}{\left(a\right)}\right)^-+c\right)\\
                &=\left(b^-\circ a^-+c\right)^-\circ \left(\left(\rhoo{b}{\left(a\right)}\right)^-+c\right)&\quad\mbox{by \eqref{eq:lambdacircrho}}\\
                &=\left(b^-\circ a^-+c\right)^-\circ \left(\left(\rhoo{b}{\left(a\right)}\right)^-+0+c\right)&\quad\mbox{by \eqref{eq:lemma24}}\\
                &=\left(b^-\circ a^-+c\right)^-\circ \left(\left(\rhoo{b}{\left(a\right)}\right)^-+\lambdaa{0}{\left(c\right)}\right)\\
                &=\left(b^-\circ a^-+\lambdaa{b^-}{\lambdaa{b}{\left(c\right)}}\right)^-\circ \left(b^-\circ\left(a^-+b\right)+\lambdaa{b^-}{\lambdaa{b}{\left(c\right)}}\right)&\quad\mbox{by \eqref{eq:lambdahomomorphism}}\\
                &=\left(a^-+\lambdaa{b}{\left(c\right)}\right)^-\circ b\circ b^-\circ\left(a^-+b+\lambdaa{b}{\left(c\right)}\right)&\quad\mbox{by \eqref{eq:semibrace}}\\
                &=\left(a^-+\lambdaa{b}{\left(c\right)}\right)^- \circ\left(a^-+b\circ\left(0+c\right)\right),&\quad\mbox{by \eqref{eq:semibrace}}
            \end{align*}
            i.e., \eqref{eq:sollambdarho} holds. Finally
            \begin{align*}
                &\rhoo{\rhoo{c}{\left(b\right)}}{\rhoo{\lambdaa{b}{\left(c\right)}}{\left(a\right)}}
                =\left(\left(\rhoo{\lambdaa{b}{\left(c\right)}}{\left(a\right)}\right)^-+\rhoo{c}{\left(b\right)}\right)^-\circ\rhoo{c}{\left(b\right)}\\
                &=\left(\left(\lambdaa{b}{\left(c\right)}\right)^-\circ \left(a^-+\lambdaa{b}{\left(c\right)}\right)+\rhoo{c}{\left(b\right)}\right)^-\circ \rhoo{c}{\left(b\right)}\\
                &=\left(\left(\lambdaa{b}{\left(c\right)}\right)^-\circ \left(a^-+\lambdaa{b}{\left(c\right)}\right)+\lambdaa{0}{\rhoo{c}{\left(b\right)}}\right)^-\circ \rhoo{c}{\left(b\right)}&\mbox{by \eqref{eq:lemma24}}\\
                &=\left(\left(\lambdaa{b}{\left(c\right)}\right)^-\circ \left(a^-+\lambdaa{b}{\left(c\right)}\right)+\lambdaa{\left(\lambdaa{b}{\left(c\right)}\right)^-}{\lambdaa{\lambdaa{b}{\left(c\right)}}{\rhoo{c}{\left(b\right)}}}\right)^-\circ \rhoo{c}{\left(b\right)}&\mbox{by \eqref{eq:lambdahomomorphism}}\\
                &=\left(\left(\lambdaa{b}{\left(c\right)}\right)^-\circ \left(a^-+\lambdaa{b}{\left(c\right)}+\lambdaa{\lambdaa{b}{\left(c\right)}}{\rhoo{c}{\left(b\right)}}\right)\right)^-\circ \rhoo{c}{\left(b\right)}&\mbox{by \eqref{eq:semibrace}}\\
                &=\left(a^- + \lambdaa{b}{\left(b\right)}\circ\left(0+\rhoo{c}{\left(b\right)}\right)\right)^-\circ\lambdaa{b}{\left(c\right)}\circ\rhoo{c}{\left(b\right)}&\mbox{by \eqref{eq:semibrace}}\\
                &=\left(a^-+b\circ\left(0+c\right)\right)^-b\circ c&\mbox{by \eqref{eq:condsolution} and \eqref{eq:lambdacircrho}}
            \end{align*}
            and
            \begin{align*}
                \rhoo{c}{\rhoo{b}{\left(a\right)}}
                &=\left(\left(\rhoo{b}{\left(a\right)}\right)^-+c\right)^-\circ c\\
                &=\left(b^-\circ\left(a^-+b\right)+c\right)^-\circ c\\
                &=\left(b^-\circ\left(a^-+b\right)+0+c\right)^-\circ c&\quad\mbox{by \eqref{eq:lemma24}}\\
                &=\left(b^-\circ\left(a^-+b\right)+\lambdaa{b^-}{\lambdaa{b}{\left(c\right)}}\right)&\quad\mbox{by \eqref{eq:lambdahomomorphism}}\\
                &=\left(b^-\circ\left(a^-+b+\lambdaa{b}{\left(c\right)}\right)\right)^-\circ c&\quad\mbox{by \eqref{eq:semibrace}}\\
                &=\left(b^-\circ\left(a^-+b\circ\left(0+c\right)\right)\right)^-\circ c,&\quad\mbox{by \eqref{eq:semibrace}}
            \end{align*}
            i.e., \eqref{eq:solrho} holds. Hence $r_B$ is a solution. Conversely if $r_B$ is a solution, in particular \eqref{eq:sollambdarho} holds. Therefore
            \begin{align*}
                &\left(a^-+\lambdaa{b}{\left(c\right)}\right)^- \circ\left(a^-+\lambdaa{b}{\left(c\right)}\circ\left(0+\rhoo{c}{\left(b\right)}\right)\right)\\&=\left(a^-+\lambdaa{b}{\left(c\right)}\right)^- \circ\left(a^-+b\circ\left(0+c\right)\right)
            \end{align*}
            holds for all $a,b,c \in B$ and since $\left(B,\circ\right)$ is a group this is equivalent to \eqref{eq:condsolution}.
        \end{proof}
    \end{theor}    
    
    \noindent As previously mentioned, this result includes Theorem 5.1 in \cite{VAJe19} that gives a sufficient, but not a necessary, condition to obtain a solution. Indeed by Proposition 2.14 in \cite{VAJe19} $\rho$ is an anti-homomorphism if and only if $c+a\circ(0+b)=c+a\circ b$ for all $a,b,c\in B$. An if $\rho$ is an anti-homomorphism then
    \begin{align*}
        a + \lambdaa{b}{\left(c\right)}\circ \left(0+\rhoo{c}{\left(b\right)}\right) = a + \lambdaa{b}{\left(c\right)}\circ\rhoo{b}{\left(c\right)} = a+b\circ c
        =a+b\circ(0+c),
    \end{align*}
    for all $a,b,c \in B$.\\
    Moreover, there exist semi-braces that do not satisfy condition \eqref{eq:condsolution}:
    an example is given by the semi-brace in \cite[Example 2.11]{VAJe19}.
    
    We now define the matched product of semi-braces, following the steps in \cite{CCoSt18x} for cancellative semi-braces.
    
    \begin{defin}
        Two semi-braces $B_1$ and $B_2$ with $\alpha: B_2 \to \Aut\left(B_1\right)$ a group homomorphism from $\left(B_2,\circ\right)$ into the automorphism group of $\left(B_1,+\right)$ and $\beta: B_1\to \Aut\left(B_2\right)$ a group homomorphism from $\left(B_1,\circ\right)$  into the automorphism group of $\left(B_2,+\right)$ such that
        \begin{align}\label{eq:mpsb1}
            \lambdaa{a}{\alphaa{}{\betaa{-1}{a}{\left(u\right)}}{}} = \alphaa{}{u}{\lambdaa{\alphaa{-1}{u}{\left(a\right)}}{}}\\\label{eq:mpsb2}
            \lambdaa{a}{\betaa{}{\alphaa{-1}{u}{\left(a\right)}}{}} = \betaa{}{a}{\lambdaa{\betaa{-1}{a}{\left(u\right)}}{}}
        \end{align}
        hold for all $a \in B_1$ and $u \in B_2$ is called a \emph{matched product system of left semi-braces}.
    \end{defin}
    
    \begin{theor}\label{th:matchedsemibrace}
        If $\left(B_1,B_2,\alpha,\beta\right)$ is a matched product system of left semi-braces, then $B_1\times B_2$ with respect to 
        \begin{align*}
            \left(a,u\right)+\left(b,v\right) &:=\left(a+b,u+v\right)\\
            \left(a,u\right)\circ\left(b,v\right) &:= \left(\alphaa{}{u}{\left(\alphaa{-1}{u}{\left(a\right)}\circ b\right)},\betaa{}{a}{\left(\betaa{-1}{a}{\left(u\right)}\circ v\right)}\right)
        \end{align*}
        is a left semi-brace called the \emph{matched product of} $B_1$ and $B_2$ (via $\alpha$ and $\beta$) and denoted by $B_1\bowtie B_2$. 
        Moreover, if $B_1$ and $B_2$ satisfy condition \eqref{eq:condsolution} then $B_1 \bowtie B_2$ satisfies the same condition and the solution associated with the matched product $B_1\bowtie B_2$ is equal to the matched product of $r_{B_1}$ and $r_{B_2}$ via $\alpha$ and $\beta$.
        \begin{proof}
            With the same proof of Theorem $9$ in \cite{CCoSt18x} it is easy to see that $B_1\bowtie B_2$ is a left semi-brace. 
            In particular $\left(B_1\times B_2,\circ\right)$ is a group with identity $\left(0, 0\right)$ and such that
            \begin{align*}
                \left(a,u\right)^- = \left(\alphaa{-1}{\bar{u}}{\left(a^-\right)}, \betaa{-1}{\bar{a}}{\left(u^-\right)}\right)
            \end{align*}
            for all $a\in B_1$ and $u \in B_2$ where for every pair $\left(a,u\right)\in B_1\times B_2$ we denote with $\bar{a}=\alphaa{-1}{u}{\left(a\right)}$ and $\bar{u}:=\betaa{-1}{a}{\left(u\right)}$. 
            Let $\left(a,u\right),\left(b,v\right),\left(c,w\right) \in B_1\times B_2$.
            Hence
            \begin{align}\label{eq:lambdasemi}
                \lambdaa{\left(a,u\right)}{\left(b,v\right)} = \left(\lambdaa{a}{\alphaa{}{\bar{u}}{\left(c\right)}}, \lambdaa{u}{\betaa{}{\bar{a}}{\left(w\right)}}\right)
            \end{align}
            \begin{align*}
            \lambdaa{\left(a,u\right)}{\left(b,v\right)}&=\left(a,u\right)\circ\left(\left(a,u\right)^-+\left(c,w\right)\right) \\
            &= \left(a,u\right)\circ\left(\alphaa{-1}{\bar{u}}{\left(a^-\right)}+c,\betaa{-1}{\bar{a}}{\left(u^-\right)}+w\right)\\
            &=\left(a\circ\alphaa{}{\bar{u}}{\left(\alphaa{-1}{\bar{u}}{\left(a^-\right)}+c\right)}, u\circ\betaa{}{\bar{a}}{\left(\betaa{-1}{\bar{a}}{\left(u^-\right)}+w\right)} \right)\\
            &=\left(a\circ \left(a^-+\alphaa{}{\bar{u}}{\left(c\right)}\right),u\circ\left(u^-+\betaa{}{\bar{a}}{\left(w\right)}\right) \right)\\
            &=\left(\lambdaa{a}{\alphaa{}{\bar{u}}{\left(c\right)}}, \lambdaa{u}{\betaa{}{\bar{a}}{\left(w\right)}}\right).
            \end{align*}
            Set $A:=\alphaa{}{u}{\lambdaa{\bar{a}}{\left(b\right)}}$ and $ U:=\betaa{}{a}{\lambdaa{\bar{u}}{\left(v\right)}}$. Then by \eqref{eq:mpsb1} the first component of $\rhoo{\left(a,u\right)}{\left(b,v\right)}=\left(\left(a,u\right)^-+\left(b,v\right)\right)^-\circ\left(b,v\right)=\left(A,U\right)^-\circ \left(a,u\right)\circ \left(b,v\right)$ is given by
            \begin{align*}
                &\bar{A}^- +     \alphaa{}{\bar{U}^-}{\lambdaa{\alphaa{-1}{\bar{U}^-}{\left(\bar{A}^-\right)}}{\left(a + \alphaa{}{u}{\lambdaa{\bar{a}}{\left(b\right)}}\right)}}\\
                &= \alphaa{-1}{\bar{U}}{\left(A^-\right)} +\alphaa{-1}{\bar{U}}{\lambdaa{\alphaa{}{\bar{U}}{\alphaa{-1}{\bar{U}}{\left(A^-\right)}}}{\left(a+\alphaa{}{u}{\lambdaa{\bar{a}}{\left(b\right)}}\right)}}
                \\&= \alphaa{-1}{\bar{U}}{\left(A^-+\lambdaa{A^-}{\left(a+\lambdaa{a}{\alphaa{}{\bar{u}}{\left(b\right)}}\right)} \right)} \\
                &= \alphaa{-1}{\bar{U}}{\left(A^-\circ a\circ \alphaa{}{\bar{u}}{\left(b\right)}\right)} \\
                &=\alphaa{-1}{\bar{U}}{\left(\left(\lambdaa{a}{\alphaa{}{\bar{u}}{\left(b\right)}}\right)^-\circ a\circ \alphaa{}{\bar{u}}{\left(b\right)} \right)} 
                \\&=
                \alphaa{-1}{\bar{U}}{\rhoo{\alphaa{}{\bar{u}}{\left(b\right)}}{\left(a\right)}} \\
                &= \alphaa{-1}{\betaa{-1}{A}{\left(U\right)}}{\rhoo{\alphaa{}{\bar{u}}{\left(b\right)}}{\left(a\right)}}
            \end{align*}
            and with the same computation the second component is 
            $\betaa{-1}{\alphaa{-1}{U}{\left(A\right)}}{\rhoo{\betaa{}{\bar{a}}{\left(v\right)}}{\left(u\right)}}$. Therefore
            \begin{align*}
                &\left(c,w\right) + \lambdaa{\left(a,u\right)}{\left(b,v\right)} \circ \left(\left(0,0\right) +\rhoo{\left(b,v\right)}{\left(a,u\right)}\right)\\
                &=\left(c,w\right) + \left(A,U\right)\circ\left(0+\alphaa{-1}{\betaa{-1}{A}{\left(U\right)}}{\rhoo{\alphaa{}{\bar{u}}{\left(b\right)}}{\left(a\right)}},0+\betaa{-1}{\alphaa{-1}{U}{\left(A\right)}}{\rhoo{\betaa{}{\bar{a}}{\left(v\right)}}{\left(u\right)}}\right)\\
                &=\left(c,w\right) +\left( \alphaa{}{U}{\left(\alphaa{-1}{U}{\left(A\right)}\circ\left(0+\alphaa{-1}{\betaa{-1}{A}{\left(U\right)}}{\rhoo{\alphaa{}{\bar{u}}{\left(b\right)}}{\left(a\right)}}\right)\right)},\right.\\&\quad\left. \betaa{}{A}{\left(\betaa{-1}{A}{\left(U\right)}\circ\left(0+\betaa{-1}{\alphaa{-1}{U}{\left(A\right)}}{\rhoo{\betaa{}{\bar{a}}{\left(v\right)}}{\left(u\right)}}\right) \right)} \right)
             \end{align*}
             hence the first component of $\left(c,w\right) + \lambdaa{\left(a,u\right)}{\left(b,v\right)} \circ \left(\left(0,0\right) +\rhoo{\left(b,v\right)}{\left(a,u\right)}\right)$ is given by
             \begin{align*}
                 &c + \alphaa{}{U}{\left(\alphaa{-1}{U}{\left(A\right)} + \lambdaa{\alphaa{-1}{U}{\left(A\right)}}{\alphaa{-1}{\betaa{-1}{A}{\left(U\right)}}{\rhoo{\alphaa{}{\bar{u}}{\left(b\right)}}{\left(a\right)}}}\right)}\\
                 &=c + A + \alphaa{}{U}{\lambdaa{\alphaa{-1}{U}{\left(A\right)}}{\alphaa{-1}{\betaa{-1}{A}{\left(U\right)}}{\rhoo{\alphaa{}{\bar{u}}{\left(b\right)}}{\left(a\right)}}}}\\
                 &=c + A + \lambdaa{A}{\rhoo{\alphaa{}{\bar{u}}{\left(b\right)}}{\left(a\right)}}\\
                 &=c + A \circ \left(0+\rhoo{\alphaa{}{\bar{u}}{\left(b\right)}}{\left(a\right)}\right)
                 \\&= c + a \circ \left(0+\alphaa{}{\bar{u}}{\left(a\right)}\right) &\qquad \mbox{by \eqref{eq:condsolution} of $B_1$}
             \end{align*}
             and with same computation the second component is $w + u \circ \left(0+\betaa{}{\bar{a}}{\left(v\right)}\right)$.\\
             On the other side
             \begin{align*}
                &\left(c,w\right) + \left(a,u\right)\circ \left(\left(0,0\right) +\left(b,v\right)\right)\\
                &=\left(c,w\right) + \left(a,u\right) \circ \left(0+b,0+v\right)\\
                &=\left(c,w\right) + \left(\alphaa{}{u}{\left(\alphaa{-1}{u}{\left(b\right)}\circ\left(0+c\right)\right)},\betaa{}{a}{\left(\betaa{-1}{a}{\left(v\right)}\circ\left(0+v\right)\right)}\right)\\
                &=\left(c,w\right) + \left(\alphaa{}{u}{\left(\alphaa{-1}{u}{\left(b\right)}+\lambdaa{\alphaa{-1}{u}{\left(b\right)}}{\left(b\right)}\right)},\betaa{}{a}{\left(\betaa{-1}{a}{\left(v\right)}+\lambdaa{\betaa{-1}{a}{\left(v\right)}}{\left(v\right)}\right)}\right)\\
                &=\left(c,w\right) + \left(b + \alphaa{}{u}{\lambdaa{\alphaa{-1}{u}{\left(b\right)}}{\left(b\right)}}, v + \alphaa{-1}{a}{\lambdaa{\betaa{-1}{a}{\left(v\right)}}{\left(v\right)}}\right)\\
                &=\left(c,w\right) + \left(b + \lambdaa{b}{\alphaa{}{\bar{u}}{\left(b\right)}}, v + \lambdaa{v}{\betaa{}{\bar{a}}{\left(v\right)}}\right)\\
                &= \left(c,w\right) + \left(b\circ\left(0+\alphaa{}{\bar{u}}{\left(b\right)}\right), v \circ \left(0 + \betaa{}{\bar{a}}{\left(v\right)}\right)\right),
            \end{align*}
            i.e., condition \eqref{eq:condsolution} holds for the semi-brace $B_1 \bowtie B_2$.
            Finally, it is clear that the solution associated to $B_1\bowtie B_2$ is actually the matched product of solutions $r_{B_1}$ and $r_{B_2}$ via $\alpha$ and $\beta$
        \end{proof}
    \end{theor}
    
    Finally, in \cite{VAJe19} the matched product of left semi-braces is defined as a generalization of the matched product of left cancellative semi-braces originally introduced in the thesis of Colazzo \cite[Theorem 3.1.1]{Cphdthesis}. In the following we show that this definition coincides with the one given in \cref{th:matchedsemibrace}.
    \begin{defin}[Definition $3.1$ in \cite{VAJe19}]
        Let $\left(B_1,+,\circ\right)$ and $\left(B_2,+,\circ\right) $ be left semi-braces. Let $\delta:B_1\to \Aut\left(B_2\right)$ be a right action of the group $\left(B_1,\circ\right)$ on the set $B_2$ and $\sigma: B_2 \to \Sym\left(B_1\right)$ a left action of the group $\left(B_2,\circ\right)$ on the set $B_1$. Assume the following properties hold for any $x, y \in B_2$ and $a, b \in B_1$:
        \begin{enumerate}
            \item $\sigma_x\left(a\circ b\right) = \sigma_x\left(a\right) \circ \sigma_{\delta_a\left(x\right)}\left(b\right)$,
            \item $\sigma_x\left(0\right) = 0$,
            \item $\delta_a\left(x\circ y\right) = \delta_{\sigma_y\left(a\right)}\left(x\right) \circ \delta_a\left(y\right)$,
            \item $\delta_a\left(0\right) = 0$,
            \item $\left(\delta_a\left(\left(x+y\right)^-\right)\right)^- = \left(\delta_a\left(x^-\right)\right)^- + \left(\delta_a\left(y^-\right)\right)^-$.
        \end{enumerate}
        Then the following operations define a left semi-brace structure on $B_1 \times B_2$
        \begin{align*}
            \left(a,x\right)+\left(b,y\right) &:=\left(a+b,x+y\right)\\
            \left(a,x\right)+\left(b,y\right) &:=\left(a\circ \sigma_{\left(\delta_{a}\left(x^-\right)\right)^-}\left(b\right), x \circ\left(\delta_{\left(\sigma_{x^-}\left(a\right)\right)^-}\right)^-\left(y^-\right) \right)
        \end{align*}
        This is called the matched product of the left semi-brace $B_1$ and $B_2$ by $\delta$ and $\sigma$.
    \end{defin}
    
    It is easy to see that the two definitions of matched product of left semi-braces coincide. Indeed, if $\left(B_1,B_2,\alpha, \beta\right)$ is a matched product system of left semi-brace it is sufficient to define $\sigma_u = \alpha_u$ and $\delta_a\left(u\right) = \betaa{-1}{\alphaa{}{u}{\left(a\right)}}{\left(u\right)}$ for all $a \in B_1$ and $u \in B_2$. Vice versa, if $B_1$, $B_2$, $\sigma$, and $\delta$ satisfy the previous definition then it is sufficient to set $\alpha_u = \sigma_u$ and $\betaa{}{a}{\left(u\right)} = \delta^{-1}_{\sigma^{-1}_u\left(a\right)}\left(u\right)$, for all $a \in B_1$ and $u \in B_2$.
	
\section{The matched product of the solutions of finite order}
	
	In this section, we focus on analyzing the matched product of solutions of finite order.
	The following lemma is a key tool to prove all results presented in this section; its proof is technical and is given at the end of the section. 
	\begin{lemma}\label{le:match-lj}
	    Let $\left(r_{S}, r_{T}, \alpha, \beta\right)$ be a matched product system of solutions. If $l\in\mathbb{N}$ and $j\in\mathbf{N}_{0}$, then $r_{S}^{l} = r_{S}^{j}$ and $r^{l}_{T} = r_{T}^{j}$ if and only if 
		$\left(r_{S}\bowtie r_{T}\right)^{l} = \left(r_{S}\bowtie r_{T}\right)^{j}$.
	\end{lemma}
	
	At first sight, \cref{le:match-lj} might seem to have restrictive assumptions. However, this is not the case, as it leads to the following powerful result.
	
	\begin{theor}\label{th:rSxrT-finord}
		Let $\left(r_{S}, r_{T}, \alpha, \beta\right)$ be a matched product system of solutions.
		Then, the solutions $r_{S}$ and $r_{T}$ are of finite order if and only if the solution $r_{S}\bowtie r_{T}$ is of finite order.
	\begin{proof}
	    First assume that $r_{S}$ and $r_{T}$ are solutions of finite order. Thus, $r_{S}^{l} = r_{S}^{k}$ and $r_{T}^{m} = r_{T}^{j}$ for certain $l,m\in\mathbb{N}$ and $i,j\in\mathbb{N}_{0}$ such that $l > k$ and $m > j$. 
	    Set $i:= mk + lj$ and $n:= lm + kj$,
	    note that $n = \left(l-k\right)\left(m-j\right) + i$
	    and $-k\left(m-j\right) + i = kj + lj\geq 0$, hence we obtain
	    \begin{align*}
	        r_{S}^{n} 
	        = r_{S}^{l\left(m-j\right)}r_{S}^{-k\left(m-j\right) + i}
	        = r_{S}^{k\left(m-j\right)}r_{S}^{-k\left(m-j\right) + i}
	        = r_{S}^{i}.
	    \end{align*}
	    Similarly, noting that $-j\left(l-k\right) + i = jk + mk\geq 0$, 
	    we obtain that $r_{T}^{n} = r_{T}^{i}$. 
	    Therefore, by \cref{le:match-lj} it holds that $\left(r_{S}\bowtie r_{T}\right)^{n} = \left(r_{S}\bowtie r_{T}\right)^{i}$ and consequently $r_{S}\bowtie r_{T}$ is a solution of finite order.
	    
	    Conversely, if $r_{S}\bowtie r_{T}$ is a solution of finite order then $\left(r_{S}\bowtie r_{T}\right)^{n} = \left(r_{S}\bowtie r_{T}\right)^{i}$ for certain $n\in\mathbb{N}$ and $i\in\mathbb{N}_{0}$ with $n > i$. By \cref{le:match-lj} it follows that $r_{S}^{n} = r_{S}^{i}$ and $r_{T}^{n} = r_{T}^{i}$ and hence both $r_{S}$ and $r_T$ are solutions of finite order.
	\end{proof}
	\end{theor}
    
    Determining the order of the matched product of two solutions of finite order requires the notion of index and period. We recall that the \emph{index} and the \emph{period} of any solution $r$ of finite order are defined as 
    \begin{align*}
        \indd{\left(r\right)}:=\min\left\{\left.j\right|\ j \in \mathbb{N}_0, \exists l \in \mathbb{N}\ r^l=r^j\right\},
    \end{align*}
    \begin{align*}
        \perr{\left(r\right)}:=\min\left\{\left.k\right|\ k \in \mathbb{N}, r^{r+\indd{\left(r\right)}}=r^{\indd{\left(r\right)}}\right\}.
    \end{align*}
    We note that
    if $h\in\mathbb{N}_{0}$, $h\geq \indd{\left(r\right)}$, then, for every $q\in\mathbb{N}_{0}$, $r^{\, \perr{\left(r\right)}q + h} = r^{h}$. 
	In addition, for all $m\in\mathbb{N}$ and $j\in\mathbb{N}_{0}$, $m > j$, it holds that $r^{m} = r^{j}$ if and only if $\perr{\left(r\right)}\mid m-j$.
	
	The following proposition allows for establishing the index and the period of the matched product of two solutions.
	\begin{prop}\label{th:rSxrT-li-mj}
		Let $\left(r_{S}, r_{T}, \alpha, \beta\right)$ be a matched product system of solutions.
		If $r_{S}$ and $r_{T}$ are solutions of finite order, then 
		\begin{align*}
		    \indd{\left(r_{S}\bowtie r_{T}\right)} = i 
		    \quad \text{and} \quad
		    \perr{\left(r_{S}\bowtie r_{T}\right)} = n - i
		\end{align*}
 		where $i:= \max\left\lbrace \indd{\left(r_{S}\right)}, \indd{\left(r_{T}\right)}\right\rbrace$ and $n:= \lcm\left(\perr{\left(r_{S}\right)}, \perr{\left(r_{T}\right)}\right) + i$.
		\begin{proof}
	        Let $q_{1}, q_{2}\in\mathbb{N}$ such that $n = \perr{\left(r_{S}\right)}q_{1} + i$ and $n = \perr{\left(r_{S}\right)}q_{2} + i$. Thus, 
	        it holds that $r_{S}^{n} = r_{S}^{i}$ and $r_{T}^{n} = r_{T}^{i}$ and by \cref{le:match-lj} we obtain that  $\left(r_{S}\bowtie r_{T}\right)^{n} = \left(r_{S}\bowtie r_{T}\right)^{i}$. 
		    Moreover, assuming that $i = \indd{\left(r_{S}\right)}$, if $\left(r_{S}\bowtie r_{T}\right)^{n} = \left(r_{S}\bowtie r_{T}\right)^{h}$ for a certain $h\in\mathbb{N}_{0}$, in particular one has that $r_{S}^{n} = r_{S}^{h}$.
		    It follows that 
		    \begin{align*}
		        r_{S}^{ \ \perr{\left(r_{S}\right)} + i} = r_{S}^{i} = r_{S}^{ \ \perr{\left(r_{S}\right)}q_{1} + i} = r_{S}^{n} = r_{S}^{h},
		    \end{align*}
		    hence $i\leq h$ and so $\indd{\left(r_{S}\bowtie r_{T}\right)} = i$. Clearly, this check is similar if one assumes that $i = \indd{\left(r_{T}\right)}$. In addition, if $\left(r_{S}\bowtie r_{T}\right)^{m} = \left(r_{S}\bowtie r_{T}\right)^{i}$ for a certain $m\in\mathbb{N}$, 
		    then $r_{S}^{m} = r_{S}^{i}$ and $r_{T}^{m} = r_{T}^{i}$. Consequently, $\perr{\left(r_{S}\right)}\mid m - i$ and $\perr{\left(r_{T}\right)}\mid m - i$, thus 
		    $\lcm\left(\perr{\left(r_{S}\right)}, \perr{\left(r_{T}\right)}\right)\mid m - i$, i.e., $n - i\mid m - i$. Therefore $n-i\leq m-i$ and hence $\perr{\left(r_{S}\bowtie r_{T}\right)} = n - i$.
	    \end{proof}
	\end{prop}
    
    The index and the period of the matched product solution $r_{S}\bowtie r_{T}$ give us upper bounds of the indexes ad periods of $r_S$ and $r_T$.
    Indeed, assuming $i:= \indd{\left(r_{S}\bowtie r_{T}\right)}$ and $p:= \perr{\left(r_{S}\bowtie r_{T}\right)}$, \cref{le:match-lj} implies that $r_{S}^{p + i} = r_{S}^{i}$ and $r_{T}^{p + i} = r_{T}^{i}$. Therefore, $r_{S}$ and $r_{T}$ are both solutions of finite order. 
    Clearly, $\indd{\left(r_{S}\right)}$ and  $\indd{\left(r_{T}\right)}$ are less than $i$, and $\perr{\left(r_{S}\right)}$ and $\perr{\left(r_{T}\right)}$ divide $p$.
	\begin{prop}\label{th:rSxrT-li-mj-2}
        Let $\left(r_{S}, r_{T}, \alpha, \beta\right)$ be a matched product system of solutions.
        If $r_{S}\bowtie r_{T}$ is a solution of finite order, 
        then it holds that
        \begin{align*}
        &\perr{\left(r_{S}\right)}\mid \perr{\left(r_{S}\bowtie r_{T}\right)}
        \quad \text{and} \quad
        \perr{\left(r_{T}\right)}\mid \perr{\left(r_{S}\bowtie r_{T}\right)}\\
        &\indd{\left(r_{S}\right)}\leq \indd{\left(r_{S}\bowtie r_{T}\right)}
        \quad \text{and} \quad
        \indd{\left(r_{S}\right)}\leq \indd{\left(r_{S}\bowtie r_{T}\right)}.
        \end{align*}
    \end{prop}
    
    The following corollary is a direct consequence of \cref{th:rSxrT-li-mj} and \cref{th:rSxrT-li-mj-2}. We note that this result includes the particular case of involutive solutions and the one of idempotent solutions already considered in \cite[Corollary 5]{CCoSt18x}.
	\begin{cor}\label{cor:rSxrT-lm2}
		Let $\left(r_{S}, r_{T}, \alpha, \beta\right)$ be a matched product system of solutions. Then the following hold:
		\begin{enumerate}
		    \item $r_{S}^{l} = \id$ and $r_{T}^{m} = \id$, for certain $l,m\in\mathbb{N}$, if and only if \  
		    $\left(r_{S}\bowtie r_{T}\right)^{n} = \id$, for a certain $n\in\mathbb{N}$;
		    \item $r_{S}^{l} = r_{S}$ and $r_{T}^{m} = r_{T}$, for certain $l,m\in\mathbb{N}$, if and only if \ $\left(r_{S}\bowtie r_{T}\right)^{n} = r_{S}\bowtie r_{T}$ for a certain $n\in\mathbb{N}$.
		\end{enumerate}
    \end{cor}
    
    The following corollary shows that, under mild assumptions, the solution associated to a semi-brace has index $1$.
    In particular, our result improves \cite[Theorem 3.2]{VAJe19}.
    \begin{cor}\label{cor:fcs1}
	    Let $B$ be a completely simple left semi-brace such that $\rho$ is an anti-homomorphism and $r_{G}$ the solution associated to the skew left brace \ $G = 0 + B + 0$.
	    Thus, for every $n\in\mathbb{N}$
	    \begin{align*}
	    r_{G}^{n} = \id 
	    \ \Longleftrightarrow \ 
	    r_{B}^{n+1} = r_{B}.
	    \end{align*}
	   In particular, $G$ is a left brace if and only if \ $r_{B}^{3} = r_{B}$.
	   \begin{proof}
	   At first note that by \cite[Theorem 3.2]{VAJe19} the left semi-brace $B$ can be written as the matched product 
	    \begin{align*}
	        B = F\bowtie \left(G\bowtie E\right),
	    \end{align*}
	    where $F$ is a left semi-brace with additive structure a left zero-semigroup, $G$ is a skew left brace, and $E$ is a left semi-brace with additive structure a right zero-semigroup.
	    In addition, by \cref{th:matchedsemibrace} it holds that 
	    \begin{align*}
	        r_{B} = r_{F}\bowtie\left(r_{G}\bowtie r_{E}\right)
	    \end{align*}
	    where in particular the solution $r_{G}$ associated to $G$ is bijective (see \cite[Theorem 3.1]{GVe17}), the solution $r_{E}$ associated to $E$ is idempotent, and the solution $r_{F}$ associated to $F$ is idempotent.
	    Consequently, assuming $r_{G}^{n} = \id$, we have that $n = \perr{\left(r_{G}\right)}q$ for a certain $q\in\mathbb{N}$ and by \cref{th:rSxrT-li-mj} we obtain that $r_{B}^{n + 1} = r_{B}$. 
	    Conversely, if $r_{B}^{n + 1} = r_{B}$, by \cref{le:match-lj} we have in particular that $r_{G}^{n+1} = r_{G}$ and by the bijectivity of $r_{G}$ clearly it follows $r_{G}^{n} = \id$.
	    
	   In particular, note that $G$ is a left brace if and only if $r_{G}$ is involutive and so by what we have just proved we obtain that $G$ is a left brace if and only if $r_{B}^{3} = r_{B}$.
	   \end{proof}
    \end{cor}
    
    The question arises whether it is feasible to find solutions with index greater than $1$.
    The answer is yes: one can consider the Lyubashenko's solution \cite{Dr90}.
    Indeed, if $f$ is a map of index $\indd{} > 1$ and period $\perr{}$ from a set $X$ into itself (see for instance \cite[p. 12]{Ho95}) and $\tau$ is the twist map on $X \times X$, then the Lyubashenko's solution $r:X\times X\to X\times X$ defined by $r:= \tau\left(f\times f\right)$ 
    is of finite order. 
    Since the maps $\tau$ and  $f\times f$ commute, it holds that $r^{k} = \tau^{k}\left(f\times f\right)^{k}$, for every $k\in\mathbb{N}_{0}$.  
    Hence, if $\perr{}$ is even then $r^{\perr{} + \indd{}} = r^{\ind{}}$ and in particular \ $\indd{\left(r\right)} = \indd{}$ \ and \ $\perr{\left(r\right)} = \per{}$. 
    In the event that $\perr{}$ then $r^{2\perr{} + \indd{}} = r^{\indd{}}$ and in particular \ $\indd{\left(r\right)} = \indd{}$ \ and \ $\perr{\left(r\right)} = 2\perr{}$.
    
    We conclude this section by presenting the complete proof of \cref{le:match-lj}.
  
	\begin{proof}[{Proof of \cref{le:match-lj}}]
	For the sake of simplicity, if $r$ is a solution on a set $X$, we denote
	\begin{align*}
		\lambda^{\left(0\right)}_{x}\left(y\right)
		:= x 
		\quad \text{and} \quad
		\rho^{\left(0\right)}_{y}\left(x\right):= y
	\end{align*}
	and, for every $n\in\mathbb{N}$,
	\begin{align*}
		\lambdai{\left(n\right)}{x}{\left(y\right)}:= \lambdaa{\lambdai{\left(n-1\right)}{x}{\left(y\right)}}{\rhoi{\left(n-1\right)}{y}{\left(x\right)}}
		\quad \text{and} \quad
		\rhoi{\left(n\right)}{y}{\left(x\right)}:=
		\rhoo{\rhoi{\left(n-1\right)}{y}{\left(x\right)}}{\lambdai{\left(n-1\right)}{x}{\left(y\right)}},
	\end{align*}
	 for all $x,y\in X$. Then, it is a routine computation to verify that for every $n\in\mathbb{N}$ and for all $x,y\in X$ it holds
	\begin{align*}
		r^{n}\left(x, y\right) = \left(\lambdai{\left(n\right)}{x}{\left(y\right)}, \  \rhoi{\left(n\right)}{y}{\left(x\right)}\right).
	\end{align*}
    In particular, note that if $n,m\in\mathbb{N}_{0}$ then $r^{n} = r^{m}$ if and only if $\lambdai{\left(n\right)}{x}{\left(y\right)} =
	\lambdai{\left(m\right)}{x}{\left(y\right)}$ and $\rhoi{\left(n\right)}{y}{\left(x\right)} = \rhoi{\left(m\right)}{y}{\left(x\right)}$, for all $x,y\in X$. 
    Furthermore, for every $n\in\mathbb{N}_{0}$, 
    $\left(r_{S}\bowtie r_{T}\right)^{n}\left(\left(a,u\right), \left(b,v\right)\right)$ 
    can be expressed in the following way
		\begin{align*}
			\left(\left(\lambdai{\left(n\right)}{a}{\left(c\right)},
			\lambdai{\left(n\right)}{u}{\left(w\right)}\right),
			\left(\alphaa{-1}{\overline{U^{\left(n\right)}}}{\rhoi{\left(n\right)}{c}{\left(a\right)}}, 
			\betaa{-1}{\overline{A^{\left(n\right)}}}{\rhoi{\left(n\right)}{w}{\left(u\right)}}\right)\right),
		\end{align*}
		where $c:= \alphaa{}{\bar{u}}{\left(b\right)}$, $w:= \betaa{}{\bar{a}}{\left(v\right)}$, $A^{\left(n\right)}:= \lambdai{\left(n\right)}{a}{\left(c\right)}$, and $U^{\left(n\right)}:= \lambdai{\left(n\right)}{u}{\left(w\right)}$. 
		We prove this by induction on $n$.
	    The case $n = 0$ follows from the fact that $\left(\lambdai{\left(0\right)}{a}{\left(c\right)},
		\lambdai{\left(0\right)}{u}{\left(w\right)}\right) = \left(a,u \right)$ and from
		\begin{align*}
			\left(\alphaa{-1}{\overline{U^{\left(0\right)}}}{\rhoi{\left(0\right)}{c}{\left(a\right)}}, 
			\betaa{-1}{\overline{A^{\left(0\right)}}}{\rhoi{\left(0\right)}{w}{\left(u\right)}}\right)
			=
			\left(\alphaa{-1}{\bar{u}}{\left(c\right)}, 
			\betaa{-1}{\bar{a}}{\left(w\right)}\right)
			= \left(b,v\right)
	    \end{align*}
		since $\overline{A^{\left(0\right)}} = \alphaa{-1}{U^{\left(0\right)}}{\left(A^{\left(0\right)}\right)} = \alphaa{-1}{u}{\left(a\right)} = \bar{a}$ and similarly $\overline{U^{\left(0\right)}} = \bar{u}$. 
		Suppose that the equality holds for $n > 0$, i.e.,
		\begin{align*}
			\lambdai{\left(n\right)}{\left(a,u\right)}{\left(b,v\right)}
			= \left(\lambdai{\left(n\right)}{a}{\left(c\right)}, \lambdai{\left(n\right)}{u}{\left(w\right)}\right)
		\end{align*}
		and
		\begin{align*} 
			\rhoi{\left(n\right)}{\left(b,v\right)}{\left(a,u\right)}
			= \left(
			\alphaa{-1}{\overline{U^{\left(n\right)}}}{\rhoi{\left(n\right)}{c}{\left(a\right)}}, 
			\betaa{-1}{\overline{A^{\left(n\right)}}}{\rhoi{\left(n\right)}{w}{\left(u\right)}}\right),
		\end{align*}
		for all $\left(a,u\right), \left(b,v\right)\in S\times T$. Thus, by induction hypothesis we have
		\begin{align*}
				\lambdai{\left(n+1\right)}{\left(a,u\right)}{\left(b,v\right)}
				&= \lambdai{}{\lambdai{\left(n\right)}{\left(a,u\right)}{\left(b,v\right)}}{\rhoi{\left(n\right)}{\left(b,v\right)}{\left(a,u\right)}}\\
				&= \lambdai{}{\left(\lambdai{\left(n\right)}{a}{\left(c\right)}, \lambdai{\left(n\right)}{u}{\left(w\right)}\right)}{\left(
					\alphaa{-1}{\overline{U^{\left(n\right)}}}{\rhoi{\left(n\right)}{c}{\left(a\right)}}, 
					\betaa{-1}{\overline{A^{\left(n\right)}}}{\rhoi{\left(n\right)}{w}{\left(u\right)}}\right)}\\
				&= \left(\lambdai{}{\lambdai{\left(n\right)}{a}{\left(c\right)}}{\left(
					\alphaa{}{\overline{U^{\left(n\right)}}}{\alphaa{-1}{\overline{U^{\left(n\right)}}}{\rhoi{\left(n\right)}{c}{\left(a\right)}}}\right)}
				,
				\lambdai{}{\lambdai{\left(n\right)}{u}{\left(w\right)}}{\betaa{}{\overline{A^{\left(n\right)}}}{\betaa{-1}{\overline{A^{\left(n\right)}}}{\rhoi{\left(n\right)}{w}{\left(u\right)}}}}
				\right)\\
				&= \left(\lambdai{}{\lambdai{\left(n\right)}{a}{\left(c\right)}}{\rhoi{\left(n\right)}{c}{\left(a\right)}}
				,
				\lambdai{}{\lambdai{\left(n\right)}{u}{\left(w\right)}}{\rhoi{\left(n\right)}{w}{\left(u\right)}}
				\right)
				= \left(\lambdai{\left(n+1\right)}{a}{\left(c\right)}, \lambdai{\left(n+1\right)}{u}{\left(w\right)}\right).
			\end{align*}
			Set 
			$\mathcal{A}:= \lambdai{}{A^{\left(n\right)}}{\alphaa{}{\overline{U^{\left(n\right)}}}{\left(\alphaa{-1}{\overline{U^{\left(n\right)}}}{\rhoi{\left(n\right)}{c}{\left(a\right)}}\right)}}$ 
			and 
			$\mathcal{U}:= \lambdai{}{U^{\left(n\right)}}{\betaa{}{\overline{A^{\left(n\right)}}}{\left(\betaa{-1}{\overline{A^{\left(n\right)}}}{\rhoi{\left(n\right)}{w}{\left(u\right)}}\right)}}$, it follows that
			\begin{align*}
				\rhoi{\left(n+1\right)}{\left(b,v\right)}{\left(a,u\right)}
				&=
				\rhoi{}{\rhoi{\left(n\right)}{\left(b,v\right)}{\left(a,u\right)}}{\lambdai{\left(n\right)}{\left(a,u\right)}{\left(b,v\right)}}\\
				&= \rhoi{}{\left(\alphaa{-1}{\overline{U^{\left(n\right)}}}{\rhoi{\left(n\right)}{c}{\left(a\right)}}, \betaa{-1}{\overline{A^{\left(n\right)}}}{\rhoi{\left(n\right)}{w}{\left(u\right)}} \right)}{\left(\lambdai{\left(n\right)}{a}{\left(c\right)}, \lambdai{\left(n\right)}{u}{\left(w\right)}\right)}\\
				&= \left(\alphaa{-1}{\overbar{\mathcal{U}}}{\rhoi{}{\alphaa{}{\overline{U^{\left(n\right)}}}{\alphaa{-1}{\overline{U^{\left(n\right)}}}{\rhoi{\left(n\right)}{c}{\left(a\right)}}}}{\lambdai{\left(n\right)}{a}{\left(c\right)}}},
				\betaa{-1}{\overline{\mathcal{A}}}{\rhoi{}{\betaa{}{\overline{A^{\left(n\right)}}}{\alphaa{-1}{\overline{A^{\left(n\right)}}}{\rhoi{\left(n\right)}{w}{\left(u\right)}}}}{\lambdai{\left(n\right)}{u}{\left(w\right)}}}
				\right)\\
				&= \left(\alphaa{-1}{\overline{U^{\left(n+1\right)}}}{\rhoi{}{\rhoi{\left(n\right)}{c}{\left(a\right)}}{\lambdai{\left(n\right)}{a}{\left(c\right)}}}, 
				\betaa{-1}{\overline{A^{\left(n+1\right)}}}{\rhoi{}{\rhoi{\left(n\right)}{w}{\left(u\right)}}{\lambdai{\left(n\right)}{u}{\left(w\right)}}}
				\right)\\
				&= \left(\alphaa{-1}{\overline{U^{\left(n+1\right)}}}{\rhoi{\left(n+1\right)}{c}{\left(a\right)}}, 
				\betaa{-1}{\overline{A^{\left(n+1\right)}}}{\rhoi{\left(n+1\right)}{w}{\left(u\right)}}\right) 
			\end{align*}
			where the second last equality holds since
			$\mathcal{A} = \lambdai{}{\lambdai{\left(n\right)}{a}{\left(c\right)}}{\rhoi{\left(n\right)}{c}{\left(a\right)}} = \lambdai{\left(n+1\right)}{a}{\left(c\right)} = A^{\left(n+1\right)}$
			and similarly 
			$\mathcal{U} = U^{\left(n+1\right)}$. 
			Therefore, the claim follows.
			\bigskip
		
		Now let us prove the statement of Lemma.
		At first suppose that $r_{S}^{l} = r_{S}^{j}$ and $r^{l}_{T} = r_{T}^{j}$. Thus, $\lambdai{\left(l\right)}{a}{\left(c\right)} = \lambdai{\left(j\right)}{a}{\left(c\right)}$ and $\rhoi{\left(l\right)}{c}{\left(a\right)} = \rhoi{\left(j\right)}{c}{\left(a\right)}$, and hence one has that
		\begin{align*}
			\lambdai{\left(l\right)}{\left(a,u\right)}{\left(b,v\right)}
			&= \left(\lambdai{\left(l\right)}{a}{\left(c\right)},
			\lambdai{\left(l\right)}{u}{\left(w\right)}\right)\\
			&= \left(\lambdai{\left(j\right)}{a}{\left(c\right)},
			\lambdai{\left(j\right)}{u}{\left(w\right)}\right)\\
			&= \lambdai{\left(j\right)}{\left(a,u\right)}{\left(b,v\right)}
		\end{align*}
		and also
		\begin{align*}
			\rhoi{\left(l\right)}{\left(b,v\right)}{\left(a,u\right)}
			&=  \left(\alphaa{-1}{\overline{U^{\left(l\right)}}}{\rhoi{\left(l\right)}{c}{\left(a\right)}}, 
			\betaa{-1}{\overline{A^{\left(l\right)}}}{\rhoi{\left(l\right)}{w}{\left(u\right)}}\right)\\
			&= \left(\alphaa{-1}{\overline{U^{\left(j\right)}}}{\rhoi{\left(j\right)}{c}{\left(a\right)}}, 
			\betaa{-1}{\overline{A^{\left(j\right)}}}{\rhoi{\left(j\right)}{w}{\left(u\right)}}\right)\\
			&= \rhoi{\left(j\right)}{\left(b,v\right)}{\left(a,u\right)},
		\end{align*}
		since $\overline{A^{\left(l\right)}} = \alphaa{-1}{U^{\left(l\right)}}{A^{\left(l\right)}} = \alphaa{-1}{U^{\left(j\right)}}{A^{\left(j\right)}} = \overline{A^{\left(j\right)}}$ and similarly $\overline{U^{\left(l\right)}} = \overline{U^{\left(j\right)}}$.
		Consequently, $\left(r_{S}\bowtie r_{T}\right)^{l} = \left(r_{S}\bowtie r_{T}\right)^{j}$.
		
		Conversely, if $\left(r_{S}\bowtie r_{T}\right)^{l} = \left(r_{S}\bowtie r_{T}\right)^{j}$, set $d:= \alphaa{-1}{\bar{u}}{\left(b\right)}$ and $z:= \betaa{-1}{\bar{a}}{\left(v\right)}$ 
		it follows that 
		\begin{align*}
			\left(\lambdai{\left(l\right)}{a}{\left(b\right)}, \lambdai{\left(l\right)}{u}{\left(v\right)}\right)
			= \lambdai{\left(l\right)}{\left(a,u\right)}{\left(d,z\right)} 
			= \lambdai{\left(j\right)}{\left(a,u\right)}{\left(d,z\right)}
			= \left(\lambdai{\left(j\right)}{a}{\left(b\right)}, \lambdai{\left(j\right)}{u}{\left(v\right)}\right)
		\end{align*}
		and so $\lambdai{\left(l\right)}{a}{\left(b\right)} = \lambdai{\left(j\right)}{a}{\left(b\right)}$ and $\lambdai{\left(l\right)}{u}{\left(v\right)} = \lambdai{\left(j\right)}{u}{\left(v\right)}$. 
		Moreover, set 
		$\mathcal{A}:= \lambdai{\left(j\right)}{a}{\alphaa{}{\bar{u}}{\left(d\right)}}$ 
		and $\mathcal{U}:= \lambdai{\left(j\right)}{u}{\betaa{}{\bar{a}}{\left(z\right)}}$ 
		it holds
		\begin{align*}
			&\left(\alphaa{-1}{\overline{U^{\left(l\right)}}}{\rhoi{\left(l\right)}{b}{\left(a\right)}},
			\betaa{-1}{\overline{A^{\left(l\right)}}}{\rhoi{\left(l\right)}{v}{\left(u\right)}}\right)
			= \rhoi{\left(l\right)}{\left(d,z\right)}{\left(a,u\right)}
			= \rhoi{\left(i\right)}{\left(d,z\right)}{\left(a,u\right)}\\
			&= \left(\alphaa{-1}{\overline{\mathcal{U}}}{\rhoi{\left(j\right)}{\alphaa{}{\bar{u}}{\left(d\right)}}{\left(a\right)}},
			\betaa{-1}{\overline{\mathcal{A}}}{\rhoi{\left(j\right)}{\betaa{}{\bar{a}}{\left(z\right)}}{\left(u\right)}}
			\right)
			= \left(\alphaa{-1}{\overline{\mathcal{U}}}{\rhoi{\left(j\right)}{b}{\left(a\right)}},
			\betaa{-1}{\overline{\mathcal{A}}}{\rhoi{\left(j\right)}{v}{\left(u\right)}}
			\right).
		\end{align*}
		Note that, $\mathcal{A} = \lambdai{\left(j\right)}{a}{\alphaa{}{\bar{u}}{\left(d\right)}} = \lambdai{\left(j\right)}{a}{\left(b\right)} = \lambdai{\left(l\right)}{a}{\left(b\right)} = A^{\left(l\right)}$ 
		and similarly 
		$\mathcal{U} = U^{\left(l\right)}$. 
		By bijectivity of $\alpha_{\overline{U^{\left(l\right)}}}$ and $\beta_{\overline{A^{\left(l\right)}}}$, one obtains that $\rhoi{\left(l\right)}{b}{\left(a\right)} = \rhoi{\left(j\right)}{b}{\left(a\right)}$ and $\rhoi{\left(l\right)}{v}{\left(u\right)} = \rhoi{\left(j\right)}{v}{\left(u\right)}$. Therefore, $r^{l}_{S} = r_{S}^{j}$ and $r^{l}_{T} = r_{T}^{j}$.
	\end{proof}
	
	%--------------------------------------
	\vspace{5mm}
	
	\bibliographystyle{elsart-num-sort}
	%--------------------------------------
	\bibliography{bibliography}

\end{document}